\documentclass{cmmse2011}
\usepackage[top=4.5cm,bottom=4cm,left=27mm,right=33mm]{geometry}
\usepackage{amssymb,graphicx}

\begin{document}

%\pagestyle{empty}

%  Include here your own macros.

\newcommand{\Eqref}[1]{(\ref{#1})}

%  Data for the headings. Please fill both fields.

\markboth
  {Multiple Prior Method for the BHDP} % abridged title
  {A. Zambelli}

%  Full title. Use '\\' to force a line break.

\title{A Multiple Prior Monte Carlo Method for\\
  the Backward Heat Diffusion Problem}

%  Full author information.

%\author{Full name}{e-mail}{affiliation number (below)}

\author{Antoine E. Zambelli}{antoine\_elabdouni@berkeley.edu}{1}

%\affiliation{number}{department}{university}

\affiliation{1}{Department of Mathematics}
  {University of California, Berkeley}

%  Abstract, key words and MSC codes.

\begin{abstract}
We consider the nonlinear inverse problem of reconstructing the heat conductivity of a cooling fin, modeled by a $2$-dimensional steady-state equation with Robin boundary conditions. The Metropolis Hastings Markov Chain Monte Carlo algorithm is studied and implemented, as well as the notion of priors. By analyzing the results using certain trial conductivities, we formulate several distinct priors to aid in obtaining the solution. These priors are associated with different identifiable parts of the reconstruction, such as areas with vanishing, constant, or varying slopes. Although more research is required for some non-constant conductivities, we believe that using several priors simultaneously could help in solving the problem.
\keywords Inverse Problems, Heat Diffusion, Monte Carlo, Prior.
\end{abstract}

%
%  Main text of the article.
%
\section{Introduction}
In this problem, we attempt to reconstruct the \emph{conductivity} $K$ in a steady state heat equation of the cooling fin on a CPU. The heat is dissipated both by conduction along the fin and by convection with the air, which gives rise to our equation (with $H$ for convection, $K$ for conductivity, $\delta$ for thickness and $u$ for temperature):

\begin{equation}\label{eq:heatpde}
u_{xx}+u_{yy}=\frac{2H}{K\delta}u
\end{equation}
The CPU is connected to the cooling fin along the bottom half of the left edge of the fin. We use the Robin Boundary Conditions (detailed in \cite{sauer}):

\begin{equation}\label{eq:robinbc}
Ku_{normal}=Hu
\end{equation}
Our data in this problem is the set of boundary points of the solution to (\ref{eq:heatpde}), which we compute using a standard finite difference scheme for an $n \times m$ mesh (here $10 \times 10$ or $20 \times 20$). We denote the correct value of $K$ by $K_{\textrm{correct}}$ and the data by $d$. In order to reconstruct $K_{\textrm{correct}}$, we will take a guess $K'$, solve the forward problem using $K'$ and compare those boundary points to $d$ by implementing the Metropolis-Hastings Markov Chain Monte Carlo algorithm (or MHMCMC). Priors will need to be established to aid in the reconstruction, as comparing the boundary points alone is insufficient.

%%%%%%%%%%%%%%%%%%%%%%%%%%%%%%%%%%%%%%%%%%%%%%%%%%%%%%%%%%%%%%%%%%%%%%%%%%%%%%%%%%%%%%%%%%%%%%%%%%%%%%%%%%%%%%%%%%%%%%%%%%%%%%%%%%%%%%%%%%%%%%%%%%%%%%%%%%%%%%%%%%%%%%%%%%%%%%%%%%%%%%%%%%%%%%%%%%%%%%%%%%%%%%%
                                                                    %Section II%
%%%%%%%%%%%%%%%%%%%%%%%%%%%%%%%%%%%%%%%%%%%%%%%%%%%%%%%%%%%%%%%%%%%%%%%%%%%%%%%%%%%%%%%%%%%%%%%%%%%%%%%%%%%%%%%%%%%%%%%%%%%%%%%%%%%%%%%%%%%%%%%%%%%%%%%%%%%%%%%%%%%%%%%%%%%%%%%%%%%%%%%%%%%%%%%%%%%%%%%%%%%%%%%%

\section{MHMCMC}
Markov Chains produce a probability distribution of possible solutions (in this case conductivities) that are most likely given the observed data (the probablility of reaching the next step in the chain is entirely determined by the current step). The algorithm is as follows (\cite{fox}). Given $K_n$, $K_{n+1}$ can be found using the following:
\begin{enumerate}
\item Generate a candidate state $K'$ from $K_n$ with some distribution $g(K'|K_n)$. We can pick any $g(K'|K_n)$ so long as it satisfies
\begin{enumerate}
\item $g(K'|K_n)=0 \Rightarrow g(K_n|K')=0$
\item $g(K'|K_n)$ is the transition matrix of Markov Chain on the state space containing $K_n,K'$.
\end{enumerate}
\item With probability \begin{equation}\label{eq:alpha} \alpha(K'|K_n) \equiv min\left\{1,\frac{Pr(K'|d)g(K_n|K')}{Pr(K_n|d)g(K'|K_n)}\right\} \end{equation} set $K_{n+1}=K'$, otherwise set $K_{n+1}=K_n$ (ie. accept or reject $K'$). Proceed to the next iteration.
\end{enumerate}
Using the probability distributions of our example, (\ref{eq:alpha}) becomes

\begin{equation} \label{eq:alpha2}
\alpha (K'|K_n) \equiv min\left\{ 1, e^{\frac{-1}{2 \sigma^2}\sum_{i,j=1}^{n,m} \left[ \left( d_{ij}-d_{ij}' \right)^2 - \left( d_{ij}-d_{n_{ij}} \right)^2  \right] }\right\} 
\end{equation}
(where $d'$ and $d_n$ denote the set of measured boundary points using $K'$ and $K_n$ respectively, and $\sigma=0.1$)\\
To simplify (\ref{eq:alpha2}), collect the constants and seperate the terms relating to $K'$ and $K_n$:

\begin{eqnarray}
&&\frac{-1}{2 \sigma^2}\sum_{i,j=1}^{n,m}{\left[ \left( d_{ij}-d_{ij}' \right)^2 - \left( d_{ij}-d_{n_{ij}} \right)^2  \right]}\\
&=& \frac{-1}{2}\sum_{i,j=1}^{n,m}{\left[ \left( \frac{d_{ij}-d_{ij}'}{\sigma} \right)^2 - \left( \frac{d_{ij}-d_{n_{ij}}}{\sigma} \right)^2  \right]}\\
&=& \frac{-1}{2}\left[ D' - D_n  \right] = f_n-f'
\end{eqnarray}
Now, (\ref{eq:alpha2}) reads

\begin{equation} \label{eq:alpha3}
\alpha (K'|K_n) \equiv min\left\{ 1, e^{ f_n - f'  }\right\}
\end{equation}
We now examine the means by which we generate a guess $K'$. If the problem consists of reconstructing a constant conductivity, we can implement a uniform change, for every iteration we take a random number $\omega$ between $-0.005$ and $0.005$ and add it to every entry in $K_n$ to obtain $K'$ (we initialize $K_0$ to a matrix of $1$s). The algorithm is highly efficient, and the reconstructed value will consistently converge to that of the solution to within $\omega$.\\
In order to approximate a nonconstant $K_{\textrm{correct}}$, the obvious choice is a pointwise change, at each iteration we add $\omega$  to a random entry of $K_n$, thus genrating $K'$. Unfortunately, systematic errors occur at the boundary points of our reconstruction (they tend to rarely change from their initial position).\\
In order to sidestep this, we use a gridwise change; change a square of the mesh (chosen at random as well) by adding $\omega$ to the four corners of said square. While this fixes the boundary problem, another major issue which arises from a non-uniform change is that the reconstruction will be marred with ``spikes", which we must iron out.

%%%%%%%%%%%%%%%%%%%%%%%%%%%%%%%%%%%%%%%%%%%%%%%%%%%%%%%%%%%%%%%%%%%%%%%%%%%%%%%%%%%%%%%%%%%%%%%%%%%%%%%%%%%%%%%%%%%%%%%%%%%%%%%%%%%%%%%%%%%%%%%%%%%%%%%%%%%%%%%%%%%%%%%%%%%%%%%%%%%%%%%%%%%%%%%%%%%%%%%%%%%%%%%%
                                                                    %Section III%
%%%%%%%%%%%%%%%%%%%%%%%%%%%%%%%%%%%%%%%%%%%%%%%%%%%%%%%%%%%%%%%%%%%%%%%%%%%%%%%%%%%%%%%%%%%%%%%%%%%%%%%%%%%%%%%%%%%%%%%%%%%%%%%%%%%%%%%%%%%%%%%%%%%%%%%%%%%%%%%%%%%%%%%%%%%%%%%%%%%%%%%%%%%%%%%%%%%%%%%%%%%%%%%%

\section{The Smoothness Prior}
To aid in ironing out the wrinkles in the reconstruction we use ``priors". Priors generally require some knowledge of the quantity we wish to find, and will add a term to (\ref{eq:alpha3}). Naturally, the more unassuming the prior, the more applicable the algorithm. This applicability will be tested as often as possible throughout these tests. The first prior compares the sum of the differences between adjacent points of $K_{\textrm{correct}}$ to those of $K'$ (keping the spikes in check), and is given by

\begin{equation}
T' =\sum_{j=1}^{n}\sum_{i=2}^{m} \left( K'(i,j)-K'(i-1,j) \right)^2 + \sum_{i=1}^{m}\sum_{j=2}^{n} \left( K'(i,j)-K'(i,j-1) \right)^2
\end{equation}
\begin{equation}
T_n =\sum_{j=1}^{n}\sum_{i=2}^{m} \left( K_n(i,j)-K_n(i-1,j) \right)^2 + \sum_{i=1}^{m}\sum_{j=2}^{n} \left( K_n(i,j)-K_n(i,j-1) \right)^2
\end{equation}
and modifying (\ref{eq:alpha3}), we obtain

\begin{equation}\label{eq:alphac}
\alpha_c (K'|K_n) \equiv min\left\{ 1, e^{f_n-f' -\lambda \left( T' - T_n \right)} \right\} 
\end{equation}
So the guess is most likely if $K_n$ and $K'$ are similarily smooth (ie, $T'\approx T_n$), if an iteration gives a $K'$ that is noticeably less smooth than the last accepted iteration, we are less likely to accept it.\\
As an initial test for the smoothness prior developed above, we attempt the gridwise change on a constant conductivity ($K_{\textrm{correct}}=1.68$, using $\lambda=100$). While we can still see the problem at the boundary points, they are limited to being a noticeable nuisance as oppposed to adamantly ruining an otherwise accurate reconstruction (whose mean comes to within $\approx 5\omega$ of $1.68$). The next step is therefore to test the algorithm on a non-constant conductivity.

%%%%%%%%%%%%%%%%%%%%%%%%%%%%%%%%%%%%%%%%%%%%%%%%%%%%%%%%%%%%%%%%%%%%%%%%%%%%%%%%%%%%%%%%%%%%%%%%%%%%%%%%%%%%%%%%%%%%%%%%%%%%%%%%%%%%%%%%%%%%%%%%%%%%%%%%%%%%%%%%%%%%%%%%%%%%%%%%%%%%%%%%%%%%%%%%%%%%%%%%%%%%%%%%
                                                                    %Section III.1%
%%%%%%%%%%%%%%%%%%%%%%%%%%%%%%%%%%%%%%%%%%%%%%%%%%%%%%%%%%%%%%%%%%%%%%%%%%%%%%%%%%%%%%%%%%%%%%%%%%%%%%%%%%%%%%%%%%%%%%%%%%%%%%%%%%%%%%%%%%%%%%%%%%%%%%%%%%%%%%%%%%%%%%%%%%%%%%%%%%%%%%%%%%%%%%%%%%%%%%%%%%%%%%%%

\subsection{Results of the Smoothness Prior on the Tilted Plane and Gaussian Well}
\noindent As a simple nonconstant trial, we look at a tilted plane with constant slope, given by
\begin{equation}\label{eq:tiltedplane10x10}
K_{\textrm{correct}}(i,j)=\frac{i+j}{20}+1;
\end{equation}
Once again, we take $K_0$ to be a matrix of all $1$s and $\lambda=100$. The boundary points again have trouble increasing from $1$ to the desired values, and in so doing lower the mean value of the reconstruction; though we still consistentily get to within about $5\omega$ of the solution (in $100000$ iterations).

\begin{figure}[h!]
  \includegraphics[height=5cm]{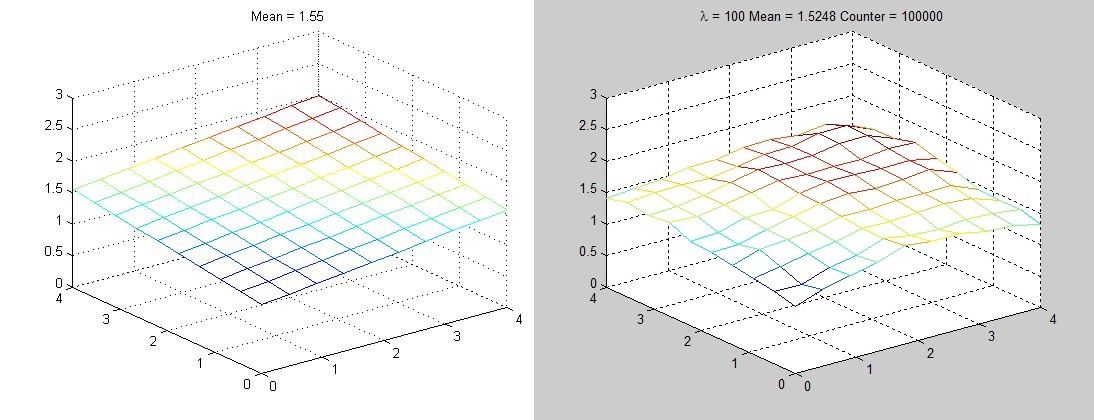}
  \caption{The $10\times10$ tilted plane we wish to reconstruct, and a reconstruction using the smoothness prior.}
\end{figure}
\noindent We now attempt to reconstruct a more complicated conductivity: a Gaussian well.\\
\newline
The Gaussian well is the first real challenge that the algorithm will face, and will be the main focus of the rest of the paper as it contains different regions which require different priors. It is given by the following equation

\begin{equation}\label{eq:gaussianwell}
K_{\textrm{correct}}(i,j) = \sum_{i=1}^{m}{\sum_{j=1}^{n}{\left(\frac{2}{1+50e^{-\frac{\left[\left( x(i) - 2 \right)^2 + \left( y(j) - 2 \right)^2 \right]}{0.2}}}\right)}}
\end{equation}
This conductivity represents a much more significant challenge, with both flat regions, and regions with steep slopes. After several trials, the optimal $\lambda$s were found to be between $1$ and $10$, though obataining a specific value for which the reconstruction is best is impossible due to the high innaccuracy of the algorithm when faced with this well.
\begin{figure}[h!]
  \includegraphics[height=5cm]{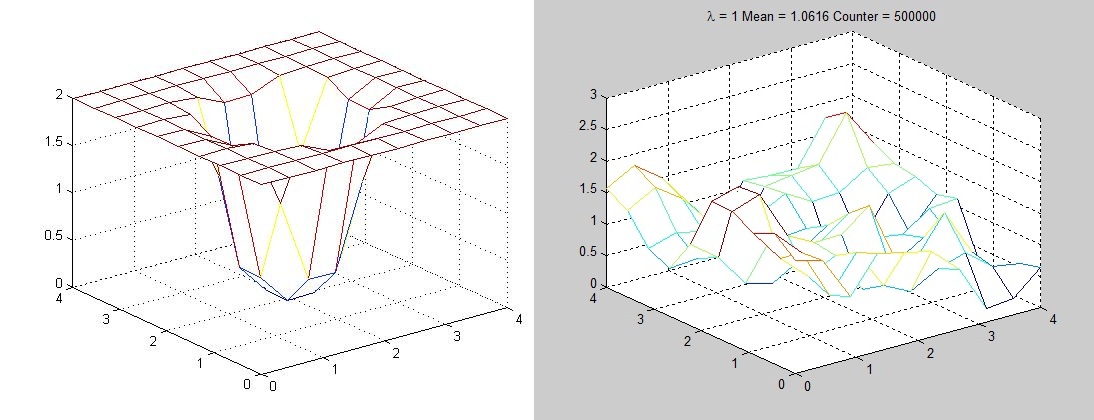}
  \caption{The $10\times10$ Gaussian well we wish to reconstruct, and a reconstruction using the smoothness prior.}
\end{figure}
\noindent There is an evident need, at this point, for much more precision. We turn once again to priors, this time developing one that will look at the slopes of the reconstruction.

%%%%%%%%%%%%%%%%%%%%%%%%%%%%%%%%%%%%%%%%%%%%%%%%%%%%%%%%%%%%%%%%%%%%%%%%%%%%%%%%%%%%%%%%%%%%%%%%%%%%%%%%%%%%%%%%%%%%%%%%%%%%%%%%%%%%%%%%%%%%%%%%%%%%%%%%%%%%%%%%%%%%%%%%%%%%%%%%%%%%%%%%%%%%%%%%%%%%%%%%%%%%%%%%
                                                                    %Section IV%
%%%%%%%%%%%%%%%%%%%%%%%%%%%%%%%%%%%%%%%%%%%%%%%%%%%%%%%%%%%%%%%%%%%%%%%%%%%%%%%%%%%%%%%%%%%%%%%%%%%%%%%%%%%%%%%%%%%%%%%%%%%%%%%%%%%%%%%%%%%%%%%%%%%%%%%%%%%%%%%%%%%%%%%%%%%%%%%%%%%%%%%%%%%%%%%%%%%%%%%%%%%%%%%%

\section{The Slope Prior}
One of the main concerns in implementing a new prior is the generality mentioned earlier. In theory, one could use a prior that only accepts Gaussian wells of the form we have here, but that code would not be very versatile. We therefore try to keep our slope prior as general as possible. In keeping with this, we look at the ratios of adjacent slopes, both in the $x$ and $y$ directions, as follows:

\begin{equation}
S_x'(i,j) = K'(i+1,j)-K'(i,j) 
\end{equation}
\begin{equation}
S_y'(i,j) = K'(i,j+1)-K'(i,j)  
\end{equation}
and define

\begin{equation}\label{eq:px}
P_x' = \sum_{j=1}^{n}{\sum_{i=1}^{m-3}{\left|\frac{S_x'(i,j)+\varepsilon_0}{S_x'(i+1,j)+\varepsilon_0}-\frac{S_x'(i+1,j)+\varepsilon_0}{S_x'(i+2,j)+\varepsilon_0}\right|}}
\end{equation}
\begin{equation}\label{eq:py}
P_y' = \sum_{j=1}^{n-3}{\sum_{i=1}^{m}{\left|\frac{S_y'(i,j)+\varepsilon_0}{S_y'(i,j+1)+\varepsilon_0}-\frac{S_y'(i,j+1)+\varepsilon_0}{S_y'(i,j+2)+\varepsilon_0}\right|}}
\end{equation}
(where $\varepsilon_0$ is $0.00005$)\\
The generality of these prior terms comes from the fact that they go to $0$ so long as the conductivity doesn't change its mind. It is equally ``happy" with a constant slope as it is with slopes that, say, double, at each grid point. It should be noted that the formulas above break down for regions where we have very small slopes adjacent to large ones, where one ratio goes to $0$ while the other grows very large. Nevertheless, we now set

\begin{equation}\label{eq:alphas}
\alpha_s (K'|K_n) \equiv min\left\{ 1, e^{f_n-f' -\mu \left( P_x'+P_y' \right)} \right\} 
\end{equation}
With this new prior, we define 

\begin{equation}
\alpha = \max\left\{\alpha_c , \alpha_s \right\}
\end{equation}
and use that in the acceptance step of the MHMCMC algorithm.

%%%%%%%%%%%%%%%%%%%%%%%%%%%%%%%%%%%%%%%%%%%%%%%%%%%%%%%%%%%%%%%%%%%%%%%%%%%%%%%%%%%%%%%%%%%%%%%%%%%%%%%%%%%%%%%%%%%%%%%%%%%%%%%%%%%%%%%%%%%%%%%%%%%%%%%%%%%%%%%%%%%%%%%%%%%%%%%%%%%%%%%%%%%%%%%%%%%%%%%%%%%%%%%%
                                                                    %Section IV.1%
%%%%%%%%%%%%%%%%%%%%%%%%%%%%%%%%%%%%%%%%%%%%%%%%%%%%%%%%%%%%%%%%%%%%%%%%%%%%%%%%%%%%%%%%%%%%%%%%%%%%%%%%%%%%%%%%%%%%%%%%%%%%%%%%%%%%%%%%%%%%%%%%%%%%%%%%%%%%%%%%%%%%%%%%%%%%%%%%%%%%%%%%%%%%%%%%%%%%%%%%%%%%%%%%

\subsection{Results of the Smoothness and Slope Priors}
\noindent Again, as a first test of the algorithm, we test it on the tilted plane.
The reconstructions reach the same precision in $100000$ iterations as we had with only the smoothness prior, so we have not yet implemented anything that is too problem-specific to the Gaussian well.\\
\newline
The initial result of the test on the well is arguably substantially better, but still rather imprecise. In an attempt to see more clearly, we make the mesh finer ($20\times 20$). In addition, we set $K_0$ to be a matrix of all $2$s. The results of the combined slope and smoothness priors are below.

\begin{figure}[h!]
  \includegraphics[height=5cm]{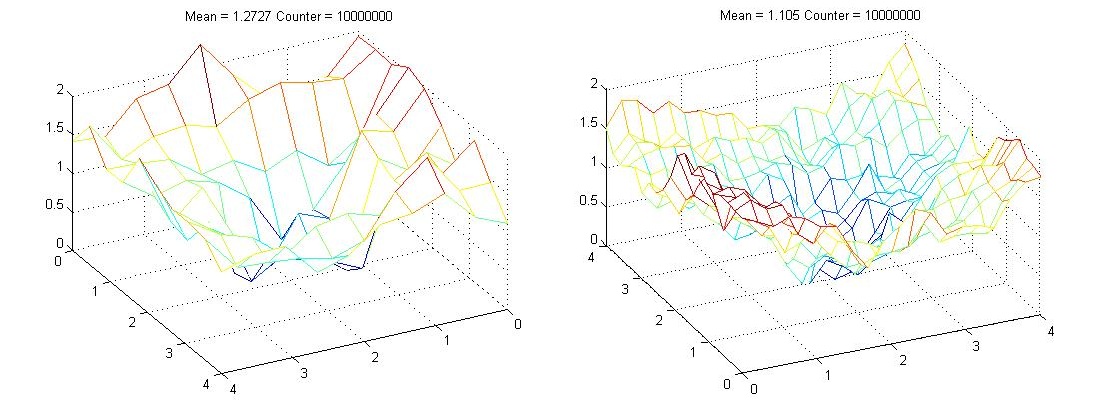}
  \caption{Results on the $10\times 10$ and $20\times20$ Gaussian wells with parameters: $\lambda=5,\ \mu=10$ and $\lambda=10,\ \mu=7.5$, respectively.}
\end{figure}
\noindent As we can see, a substantial improvement has been made over the attempt in section (3.1). We now consistently obtain somewhat of a bowl shape. In comparing the solution we wish to achieve and the reconstruction we have, one notices that the major problem areas are the outer regions, where the conductivity is nearly constant. As previously stated, equations (\ref{eq:px}) and (\ref{eq:py}) break down when the slopes are vanishing, so it is reasonable to assume that with this alone, the reconstruction will not improve as substantially as we need it to. As before, we implement another prior to aid us.

%%%%%%%%%%%%%%%%%%%%%%%%%%%%%%%%%%%%%%%%%%%%%%%%%%%%%%%%%%%%%%%%%%%%%%%%%%%%%%%%%%%%%%%%%%%%%%%%%%%%%%%%%%%%%%%%%%%%%%%%%%%%%%%%%%%%%%%%%%%%%%%%%%%%%%%%%%%%%%%%%%%%%%%%%%%%%%%%%%%%%%%%%%%%%%%%%%%%%%%%%%%%%%%%
                                                                    %Section V%
%%%%%%%%%%%%%%%%%%%%%%%%%%%%%%%%%%%%%%%%%%%%%%%%%%%%%%%%%%%%%%%%%%%%%%%%%%%%%%%%%%%%%%%%%%%%%%%%%%%%%%%%%%%%%%%%%%%%%%%%%%%%%%%%%%%%%%%%%%%%%%%%%%%%%%%%%%%%%%%%%%%%%%%%%%%%%%%%%%%%%%%%%%%%%%%%%%%%%%%%%%%%%%%%

\section{Smoothness, Flatness, and Slope Priors}
To help reconstruct the outermost regions of the well, we need a prior that will go to $0$ for regions that have vanishing slope. The most obvious choice is therefore to use what we computed for the smoothness prior

\begin{equation}
T' =\sum_{j=1}^{n}\sum_{i=2}^{m} \left( K'(i,j)-K'(i-1,j) \right)^2 + \sum_{i=1}^{m}\sum_{j=2}^{n} \left( K'(i,j)-K'(i,j-1) \right)^2 
\end{equation}
and set

\begin{equation}\label{eq:alphaf}
\alpha_f (K'|K_n) \equiv min\left\{ 1, e^{f_n-f' -W \left( T' \right)} \right\} 
\end{equation}
using

\begin{equation}
\alpha = \max\left\{\alpha_c , \alpha_s , \alpha_f \right\}
\end{equation}
in the MHMCMC algorithm. Again, the worry that adding a new prior would undermine the generality of the algorithm can be eased by noting that we are simply accounting for a problematic case not treated by the slope prior, though we still test this prior on the tilted plane.

%%%%%%%%%%%%%%%%%%%%%%%%%%%%%%%%%%%%%%%%%%%%%%%%%%%%%%%%%%%%%%%%%%%%%%%%%%%%%%%%%%%%%%%%%%%%%%%%%%%%%%%%%%%%%%%%%%%%%%%%%%%%%%%%%%%%%%%%%%%%%%%%%%%%%%%%%%%%%%%%%%%%%%%%%%%%%%%%%%%%%%%%%%%%%%%%%%%%%%%%%%%%%%%%
                                                                    %Section V.1%
%%%%%%%%%%%%%%%%%%%%%%%%%%%%%%%%%%%%%%%%%%%%%%%%%%%%%%%%%%%%%%%%%%%%%%%%%%%%%%%%%%%%%%%%%%%%%%%%%%%%%%%%%%%%%%%%%%%%%%%%%%%%%%%%%%%%%%%%%%%%%%%%%%%%%%%%%%%%%%%%%%%%%%%%%%%%%%%%%%%%%%%%%%%%%%%%%%%%%%%%%%%%%%%%

\subsection{Results of the Combined Priors}
\noindent The $20\times20$ tilted plane is given by

\begin{equation}\label{eq:tiltedplane20x20}
K_{\textrm{correct}}(i,j)=\frac{i+j}{40}+1
\end{equation}
Running the MHMCMC algorithm with all three priors yields fairly accurate reconstructions, that miss the solution by $\approx 6 \omega $. One should again note the presence of the familiar (though no less troublesome) boundary points.
\begin{figure}[h!]
  \includegraphics[height=5cm]{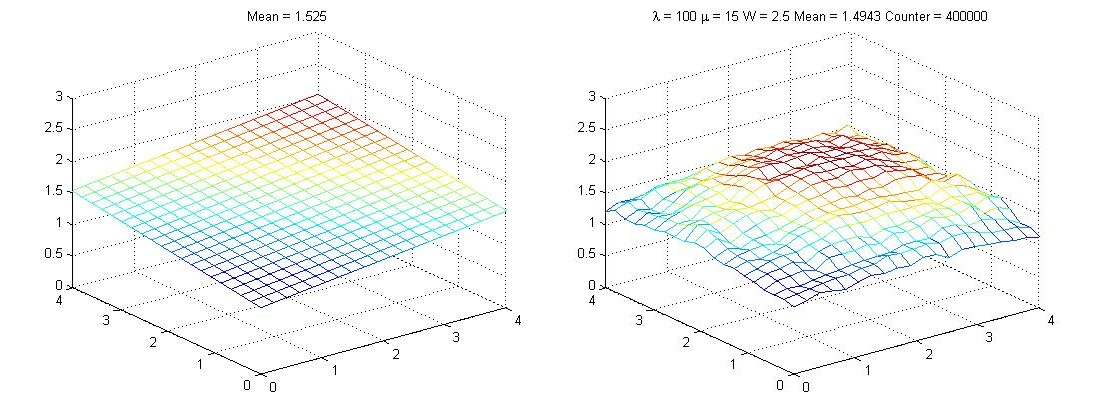}
  \caption{The $20\times20$ tilted plane, and a reconstruction using all three priors.}
  \label{fig:1}
\end{figure}
\newline
We now try once again to reconstruct the Gaussian well. The results of the added prior are apparent, and regions with vanishing slope are treated much more accurately than before. Perhaps the most successful reconstruction thus far is the following (though many more possible combinations of $\lambda$, $\mu$ and W must be explored).

\begin{figure}[h!]
  \includegraphics[height=5cm]{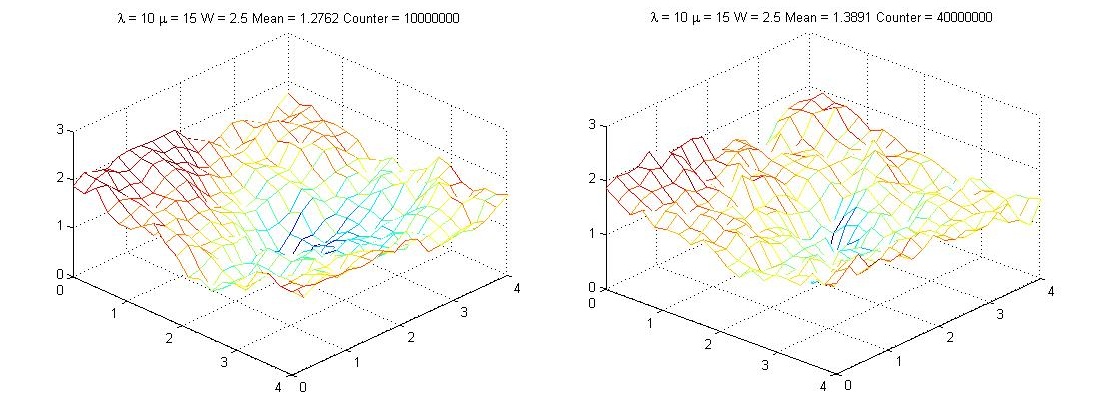}
  \caption{Results on the $20\times20$ Gaussian well using all three priors, at $10$ and $40$ million interations.}
  \label{fig:2}
\end{figure}
\noindent An obvious flaw in these reconstructions happens to be the width of the well, the algorithm is still capable of reconstructing the center of the well, and its depth, but it is often much narrower than in the actual solution. It would seem the algorithm has trouble starting to drop off from the vanishing slope region into the varying one. This exposes an inherent problem with the patchwork we have taken thus far: getting the seams to match up nicely.

%%%%%%%%%%%%%%%%%%%%%%%%%%%%%%%%%%%%%%%%%%%%%%%%%%%%%%%%%%%%%%%%%%%%%%%%%%%%%%%%%%%%%%%%%%%%%%%%%%%%%%%%%%%%%%%%%%%%%%%%%%%%%%%%%%%%%%%%%%%%%%%%%%%%%%%%%%%%%%%%%%%%%%%%%%%%%%%%%%%%%%%%%%%%%%%%%%%%%%%%%%%%%%%%
                                                                    %Section VI%
%%%%%%%%%%%%%%%%%%%%%%%%%%%%%%%%%%%%%%%%%%%%%%%%%%%%%%%%%%%%%%%%%%%%%%%%%%%%%%%%%%%%%%%%%%%%%%%%%%%%%%%%%%%%%%%%%%%%%%%%%%%%%%%%%%%%%%%%%%%%%%%%%%%%%%%%%%%%%%%%%%%%%%%%%%%%%%%%%%%%%%%%%%%%%%%%%%%%%%%%%%%%%%%%

\section{Conclusion}

As we have seen, reconstructions of the heat conductivity greatly benefit from added priors. There is certainly much work left to be done, and a very careful analysis of the seams at which the various priors trade off is in order. However, we believe that in testing the algorithm against other complex nonconstant conductivities, which is the next step we plan to take, it is possible to complete the aforementioned analysis of the seams and reconstruct complex quantities via this patchwork method.

\section*{Acknowledgements}

I was introduced to this problem at a National Science Foundation REU program at George Mason University, and would like to thank both of those institutions for the oppurtunity that gave me.
I would also like to thank Professors Timothy Sauer and Harbir Lamba at GMU, who got me started on this project while I was there and helped me decipher the MHMCMC algorithm.

%
%  Bibliography. Follow the usual conventions.
%


\begin{thebibliography}{99}

\bibitem{sauer}
  {\sc T. Sauer},
  {\em Numerical Analysis},
  Pearson Addison-Wesley, 2006.

\bibitem{fox}
  {\sc C. Fox, G.\,K. Nicholls and S.\,M. Tan},
  {\em Inverse Problems, Physics 707, The University of Auckland},
  ch 7.
  
\end{thebibliography}
\end{document}